\documentclass{article}
\usepackage{amssymb,amsfonts,latexsym}

\newtheorem{theorem}{Theorem}
\newtheorem{defi}[theorem]{Definition}
\newtheorem{remark}[theorem]{Remark}
\newtheorem{lemma}[theorem]{Lemma}
\newtheorem{prop}[theorem]{Proposition}
\newtheorem{cor}[theorem]{Corollary}

\newcommand{\e}{\emptyset}
\renewcommand{\c}{\mathfrak{C}}
\newcommand{\Aut}{\mathrm{Aut}}
\newcommand{\tp}{\mathrm{tp}}

\newcommand{\bdd}{\mathrm{bdd}}
\newcommand{\dcl}{\mathrm{dcl}}
\newcommand{\Fix}{\mathrm{Fix}}
\renewcommand{\O}{\mathcal{O}}
\newcommand{\KP}{\mathrm{KP}}
\newcommand{\sw}{\stackrel{\scriptscriptstyle \mathrm{sw}}{\sim}}

\newenvironment{proof}{\vspace{-0.25cm}
{\bf Proof}: }{\hfill $\Box$}

\begin{document}
\title{Normal hyperimaginaries\footnote{Partially supported by grants  MTM 2011-26840 of Spanish Ministry of Economy and Competitiveness and 2009SGR-187 of Catalan Government}}
\author{Enrique Casanovas and Joris Potier\\University of Barcelona}
\date{December 3, 2013}
\maketitle

\begin{abstract} We introduce the notion of normal hyperimaginary and we develop its basic theory. We present a new proof of Lascar-Pillay's theorem on bounded hyperimaginaries based on properties of normal hyperimaginaries. However, the use of Peter-Weyl's theorem on the structure of compact Hausdorff groups is not completely eliminated from the proof.  In the second part, we show that all closed sets in Kim-Pillay spaces are equivalent to hyperimaginaries and we use this to introduce an approximation of $\varphi$-types for bounded hyperimaginaries.
\end{abstract}

As usual, we work in the monster model $\c$ of a complete theory $T$ of language $L$. For background on hyperimaginaries we refer to~\cite{Casanovas09}.  Recall that a hyperimaginary is an equivalence class $e=a_E$ of a possibly infinite tuple $a$ under a $0$-type-definable equivalence relation $E$. We use the notation $E(x,y)$ for the partial type defining the equivalence relation  $E$.

For a hyperimaginary $e$, let $\Fix(e)=\Aut(\c/e)$ be the group of automorphisms of the monster model $\c$ fixing $e$. A hyperimaginary $d$ is \emph{definable over $e$} if $f(d)=d$ for all $f\in \Fix(e)$.
The \emph{definable closure} $\dcl(e)$ of $e$ is the class of all hyperimaginaries definable over $e$.
Two hyperimaginaries $e,d$ are \emph{equivalent}, written $e\sim d$, if they are interdefinable, that is, if $\dcl(e)=\dcl(d)$.  This notation can also be applied to the case where $e$ or $d$ are sequences of hyperimaginaries. If $A$ is a set of hyperimaginaries $e\sim A$ means that $e\sim d$ for a sequence $d$ enumerating $A$. In some cases we  will be interested in automorphisms fixing  $A$ set-wise.  We write $e\sw A$   to mean that $\Fix(e)$  is the set of all automorphisms $f$  such that $f(A)=A$  (set-wise).  If $(A_i:i\in I)$ is a sequence of sets, we write $e\sw (A_i:i\in I)$  meaning that  $\Fix(e)$ is the set of all automorphisms $f$ such that  $f(A)=A$    for all $i\in I$.

The \emph{cardinality} $|e|$ of a hyperimaginary $e$ is the minimal cardinality of a set $A$ of real elements (i.e., $A\subseteq \c$) such that $e\in\dcl(A)$. In this case, for any cardinal $\kappa\geq |e|$ there is a $0$-type-definable equivalence relation $E$ on $\kappa$-tuples  and there is a $\kappa$-tuple $a$ such that $e\sim a_E$. The hyperimaginary $e$ is called \emph{finitary} if $|e|<\omega$. Equivalently, $e$ is finitary if
$e\sim a_E$ for some finite tuple $a$ and some $0$-type-definable equivalence relation $E$.

A hyperimaginary $e$ is \emph{bounded} if it has a small orbit (an orbit of cardinality smaller than the size of the monster model). We denote by $\bdd(\e)$ the class of all bounded hyperimaginaries. There is a single hyperimaginary $e$ which is interdefinable with $\bdd(\e)$, in the sense that $\dcl(e)=\bdd(\e)$. More generally, for any  definably closed class $A\subseteq \bdd(\e)$ there is a single $e\in A$ such that $A=\dcl(e)$. For any index set $I$, the relation $\equiv_{\bdd(\e)}$ of having the same type over $\bdd(\e)$ restricted to $I$-tuples
is the smallest bounded (i.e., with a small number of classes) $0$-type-definable equivalence relation on $I$-tuples. It is also called the \emph{Kim-Pillay} equivalence relation and its classes are called $\KP$-strong types.  The set of all $\KP$-classes of $\alpha$-tuples is $\c^\alpha/\KP$.

We see the class of all definably closed classes  of hyperimaginaries as a lattice with the order of inclusion. Hence  $\inf(A,B)= A\cap B$ and $\sup(A,B)=\dcl(A\cup B)$. By abuse of notation we write something like $\inf(e_1,e_2)\sim d$ or even $\inf(e_1,e_2)= d$ for hypermaginaries $e_1,e_2,d$ to mean that $\inf(\dcl(e_1),\dcl(e_2))=\dcl(d)$. Note that $\sup(e_i:i\in I)=\dcl(e_i:i\in I)$.

Lascar and Pillay proved in~\cite{Las-Pil98} that every bounded hyperimaginary is equivalent to a sequence of finitary hyperimaginaries. Their proof rely on an application of Peter-Weyl's  theorem on the structure of compact Hausdorff groups according to which each such group is an inverse limit of compact Lie groups.  We seek for a purely model-theoretical proof of the same result, avoiding the use of Peter-Weyl's theorem.  There are particular cases  where the existence of such a sequence of finitary hyperimaginaries is easy to guarantee:  normal hyperimaginaries and $\KP$-classes  (see Proposition~\ref{normalfinitary} and  Lemma~\ref{l2} below)

\section{Normal hyperimaginaries}

The group $G=\Aut(\bdd(\emptyset))$ of elementary permutations of $\bdd(\emptyset)$ is a topological group, with a compact Hausdorff topology. Its closed subgroups are all subgroups of the form $\Fix_G(e)=\{f\in G: f(e)=e\}$  with $e\in \bdd(\emptyset)$. For a complete description of the topology see~\cite{Las-Pil98} or~\cite{Ziegler01}. If we endow $\Aut(\c)$ with the topology of point-wise convergence (a basis of open sets is given by all sets of the form $\{f\in\Aut(\c): f(a)= b\}$  for all finite tuples $a,b\in\c$) then $\Aut(\c)$ is a topological group and the canonical projection $\Aut(\c)\rightarrow G$ is continuous.  Notice that $G\cong \Aut(\c)/\Aut(\c/\bdd(\emptyset))$. According to Peter-Weyl's theorem, there is a family $(G_i:i\in I)$ of normal closed subgroups $G_i$ of $G$ such that  $\bigcap_{i\in I}G_i=\{1\}$ and each $G/G_i$ is a compact Lie group, and hence it has the descending chain condition (DCC) on closed subgroups.  Each $G_i$ is of the form $\Fix_G(e_i)$  for some $e_i\in \bdd(\emptyset)$.  Let $H_i =\Fix(e_i)$ be the corresponding subgroup of $\Aut(\c)$. Then $\bigcap_{i\in I}H_i =\Aut(\c/\bdd(\emptyset))$ and therefore $(e_i:i\in I)$ is interdefinable with any tuple enumerating $\bdd(\emptyset)$.  Moreover the DCC of $G/G_i$ translates as follows: there is no strictly ascending chain $(G_{i,j}:j<\omega)$ of closed subgroups $G_{i,j}\leq G_{i,j+1}$ of $G$ extending $G_i$. This explains the following definitions:

\begin{defi}\rm
A hyperimaginary $e$ is \emph{normal} if $\Fix(e)$ is a normal subgroup of $\Aut(\c)$. A hyperimaginary $e$ is \emph{DCC} if there is no sequence $(e_n:n<\omega)$ of hyperimaginaries $e_n\in \dcl(e)$ such that $e_n\in \dcl(e_{n+1})$ and $e_{n+1}\not\in \dcl(e_n)$ for each $n<\omega$.
\end{defi}

Peter-Weyl's theorem give us a sequence $(e_i:i\in I)$ of normal DCC hyperimaginaries $e_i\in\bdd(\emptyset)$ such that $(e_i:i\in I)\sim \bdd(\emptyset)$.
We will see that normal hyperimaginaries are bounded and that normal DCC hyperimaginaries are finitary. We will show that in order to prove Lascar-Pillay's theorem it is in fact enough to find a sequence $(e_i:i\in I)$ of finitary normal hyperimaginaries $e_i$ such that $(e_i:i\in I)\sim\bdd(\emptyset)$. 

\begin{defi}\rm We call \emph{Peter-Weyl's condition}  the statement that there is  a sequence $(e_i:i\in I)$ of finitary normal hyperimaginaries $e_i$ such that $(e_i:i\in I)\sim \bdd(\emptyset)$.
\end{defi}

We have not found a proof of  Peter-Weyl's condition avoiding the use of Peter-Weyl's theorem, but we can offer an easy-to-follow proof of Lascar-Pillay's theorem assuming  this condition.

\begin{prop} The following are equivalent for any hyperimaginary $e$:
\begin{enumerate}
\item $e$ is normal.
\item For any $e^\prime\equiv e$, $e^\prime \in \dcl(e)$.
\item $e\sim (f(e):f\in \Aut(\c))$
\item $e$ is equivalent to a sequence enumerating an orbit of a hyperimaginary.
\end{enumerate}
\end{prop}
\begin{proof} \emph{1} $\Leftrightarrow$ \emph{2}. By definition, $e$ is normal iff for any   $f,g\in\Aut(\c)$ such that $f(e)=e$, we have $g^{-1}fg(e)=e$, that is $f(g(e))=g(e)$. Therefore, $e$ is normal iff $\{g(e): g\in \Aut(\c)\}\subseteq \dcl(e)$.

\emph{2} $\Rightarrow$ \emph{3}. Clear, since $f(e)\equiv e$ for every $f\in\Aut(\c)$.

\emph{3} $\Rightarrow$ \emph{4}. Obvious.

\emph{4} $\Rightarrow$ \emph{2}. If $e$ is equivalent to an enumeration of an orbit and $e^\prime\equiv e$, then $e^\prime$ is equivalent to an enumeration of the same orbit and therefore $e^\prime\in\dcl(e)$.
\end{proof}

\begin{remark} Normal hyperimaginaries are bounded.
\end{remark}
\begin{proof} Let $e$ be normal. If $(e_i:i<\kappa)$ is a long enough sequence of different conjugates of $e$, then we can find $i<j<\kappa$ with $e_i\equiv_e e_j$. Since $e_i,e_j$ are definable over $e$, $e_i=e_j$, a contradiction.
\end{proof}

\begin{prop} A hyperimaginary $e$ is normal if and only if for any index set $I$, the equivalence relation $\equiv_e$ on $I$-tuples is $0$-type-definable.
\end{prop}
\begin{proof}  Let $(e_j:j\in J)$ be a (bounded) orbit equivalent to $e$.   Then $\equiv_e\, = \,\equiv_{(e_j:j\in J)}$, which is clearly invariant and type-definable, hence $0$-type-definable.

If $\equiv_e$ is $0$-type definable, then  also $\equiv_e$ as a relation between hyperimaginaries is $0$-type-definable. Let $f\in \Fix(e)$ and $g\in\Aut(\c)$ such that $g(e)=e^\prime$. Then $e^\prime\equiv_e f(e^\prime)$. If we apply $g^{-1}$ we see that $e\equiv_e g^{-1}f(e^\prime)$ and hence $g^{-1}f(e^\prime)= e$. If we apply $g$ we conclude that $f(e^\prime)=g(e)=e^\prime$. Therefore $e^\prime\in\dcl(e)$.
\end{proof}

\begin{remark} If each $e_i$ is normal, then $(e_i:i\in I)$ is normal.
\end{remark}

\begin{lemma}\label{m} Let  $e=a_E$ be normal.
\begin{enumerate}
\item  $e\sim a_{\equiv_e}$.
\item For any tuple $m$ enumerating a model, $e\sim m_{\equiv_e}$.
\end{enumerate}
\end{lemma}
\begin{proof} \emph{1}. If $e$ is normal, then $\equiv_e$ is $0$-type-definable and $a_{\equiv_e}$ is a hyperimaginary.  Assume  first $f\in\Fix(e)$. Then $a\equiv_e f(a)$ and therefore $f(a_{\equiv_e})=a_{\equiv_e}$.  For the other direction, assume now $f(a_{\equiv_e})=a_{\equiv_e}$. Then $f(a)\equiv_e a$. Since $a_E=e$, $f(a_E)=e$, that is, $f(e)=e$.

 \emph{2}. Assume $m$ enumerates a model. Clearly, $m_{\equiv_e}\in \dcl(e)$. On the other hand, if $f$ fixes $m_{\equiv_e}$ then $m\equiv_e f(m)$ and there is some $g\in \Fix(e)$ such that $g(m)=f(m)$. It follows that $fg^{-1}$ fixes point-wise a model and it is a strong automorphism, which implies it fixes every element of $\bdd(\e)$. Hence $f(e)=fg^{-1}g(e)=fg^{-1}(e)=e$.
\end{proof}

\begin{prop}\label{normalfinitary} Every normal hyperimaginary is equivalent to a sequence of finitary hyperimaginaries.
\end{prop}
\begin{proof} Let $e$ be normal. By the previous lemma, $\equiv_e$ is type-definable over $\e$ and  $e\sim a_{\equiv_e}$ for some tuple $a$. Let $a=(a_i:i<\kappa)$ and for each finite $X\subseteq \kappa$ let $E^X$ be defined  for $\kappa$-tuples $b,c$ by
$$E^X(b,c) \Leftrightarrow b\restriction X\equiv_e c\restriction X.$$
If $e^X= a_{E^X}$, then each $e^X$ is finitary and $e\sim (e^X: X\subseteq \kappa \mbox{ finite })$.
\end{proof}

\begin{lemma} Every  normal DCC hyperimaginary is finitary.
\end{lemma}
\begin{proof} Let $e$ be normal DCC.  Choose, like in the proof of Proposition~\ref{normalfinitary}, a tuple $a=(a_i:i<\kappa)$ such that $e\sim a_{\equiv_e}$ and define $E^X$ and $e^X$ as in that proof. Clearly, $e^X\in \dcl(e)$ and if $X\subseteq Y$, then $e^X\in \dcl(e^Y)$. Since $e$ is DCC, there is some finite $X$ such that for all finite $Y\supseteq X$, $e^Y\in \dcl(e^X)$. It follows that $e\sim e^X$ and hence $e$ is finitary.
\end{proof}

\begin{prop}\label{EF}
\begin{enumerate}
\item For any $0$-type-definable equivalence relation on $\kappa$-tuples $F$, for any hyperimaginary $e$, if $E=\equiv_e$, then the relational product $E\circ F=F\circ E=E\circ F\circ E$ is an equivalence relation.
\item Given normal $e$ and $d\in\bdd(\e)$, there are  a $\kappa$-tuple $m$ and a $0$-type-definable equivalence relation $F$ on $\kappa$-tuples  such that, if  $E$ is the $0$-type-definable equivalence relation $\equiv_{e}$ on $\kappa$-tuples, then  $m_E\sim e$, $m_F\sim d$ and $m_{E\circ F}\sim \inf(e,d)$.
\end{enumerate}
\end{prop}
\begin{proof} \emph{1}. We must check symmetry and transitivity of $E\circ F$. For symmetry, assume $a\equiv_e bFc$ and choose an automorphism $f$ such that $f(e)=e$ and $f(a)=b$. Let $c^\prime$ be such that $f(c^\prime)=c$. Then $ac^\prime\equiv bc$ and 
therefore $F(a,c^\prime)$. Hence $c\equiv_e c^\prime F a$. Using now  symmetry, for transitivity it is enough to prove that if  $a \equiv_e b F c \equiv_e  d$, then $a E\circ F d$.
Choose  $f\in \Fix(e)$ such that  $f(c)=d$.  Then $a\equiv_e f(b) F d$.

\emph{2}. Let $d= a_G$ for a tuple $a$, and extend  $a$  to a tuple $m=(m_i:i<\kappa)$ enumerating a model.  Let $I\subseteq \kappa$ be such that $a=(m_i:i\in I)$ and define $F$ by
$$F(x,y)\leftrightarrow G(x\restriction I,y\restriction I) .$$
It is a $0$-type-definable equivalence relation and $m_F\sim d$.  Let $E= \equiv_e$. By Lemma~\ref{m},  $m_E\sim e$.  It is  clear that $m_{E\circ F}\in \dcl(m_E)\cap \dcl(m_F)$. Now we assume $e^\prime\in \dcl(m_E)\cap \dcl(m_F)$ and we check that $e^\prime\in\dcl(m_{E\circ F})$. For this purpose, let $f$ be an automorphism fixing $m_{E\circ F}$. Then $E\circ F(m,f(m))$ and by symmetry $F\circ E(m,f(m))$. Let $b$ be such that $F(m,b)\wedge E(b,f(m))$. 
Since $b\equiv_e f(m)$, there is an automorphism $g\in\Fix(e)$ such that $g(b)=f(m)$.   Then  $F(g(m), g(b))$, that is $F(m, g^{-1}f(m))$. Let $h=g^{-1}f$. Since $h$ fixes $m_F$, $h(e^\prime)=e^\prime$. Since $g\in\Fix(e)$,  $m\equiv_e g(m)$ and hence $g$ fixes $m_E$ and $g(e^\prime)=e^\prime$. Therefore $f(e^\prime)= gh(e^\prime)= g(e^\prime) =e^\prime$.
\end{proof}

\begin{remark}Under the Galois correspondence, $\inf (e,d)$ corresponds to $\sup(\Fix(e),\Fix(d))$ in the lattice of closed subgroups. If $\Fix(e)$ is a normal subgroup, this $\sup$ is the product $\Fix(e)\cdot\Fix(d)$ (the product of two compact subgroups is compact, hence closed).
So in Proposition~\ref{EF},  $\Fix(e)\cdot\Fix(d)=\Fix(m_{E\circ F}) $. 
\end{remark}

\begin{remark}\label{W} To prove Peter-Weyl's condition it is enough to prove  that for  every finitary bounded hyperimaginary $e$ there is a family $(e_i:i\in I)$ of finitary  normal hyperimaginaries $e_i$ such that $e\in \dcl(e_i:i\in I)$.
\end{remark}
\begin{proof} There is  a normal $e$ such that $e\sim \bdd(\e)$. Since $e$ is equivalent to a family of finitary bounded hyperimaginaries and each
finitary bounded hyperimaginary is definable over a family of  finitary normal hyperimaginary, we conclude that $e$ is definable over a family $(e_i:i\in I)$ of finitary  normal hyperimaginaries. It follows that $e\sim (e_i:i\in I)$.
\end{proof}

\begin{cor}[Lascar-Pillay] \label{LP}
Every bounded hyperimaginary is equivalent to a sequence of finitary hyperimaginaries.
\end{cor}
\begin{proof}(Assuming Peter-Weyl's condition) Let $d$ be a bounded hyperimaginary and choose  a family $(e_i:i\in I)$ of finitary  normal hyperimaginaries such that $(e_i:i\in I)\sim \bdd(\e)$.  Let $\kappa \geq |I|,|d|,|T|$, and for each $i\in I$ let $E_i$ be the equivalence  relation
$\equiv_{e_i}$ on $\kappa$-tuples.  Let $E$ be the Kim-Pillay equivalence relation $\equiv_{\bdd(\e)}$ on $\kappa$-tuples. We may assume that the family is closed under finite composition (that is, for any $i,j\in I$ there is some $k\in I$ such that $e_k\sim e_ie_j$), which implies $E=\bigcap_{i\in I}E_i$. Choose with Proposition~\ref{EF} a $0$-type-definable bounded equivalence relation $F$ on $\kappa$-tuples and some $\kappa$-tuple $m$ such that $d\sim m_F$, $e_i\sim m_{E_i}$ and $\inf(e_i,d)\sim m_{E_{i}\circ F}$. Since $e_i$ is finitary, $\inf(e_i,d)$ is finitary too. We claim that $d\sim (\inf (e_i,d):i\in I)$.  Notice that  $F= E\circ F$. Hence $d\sim m_{E\circ F}$ and it is enough to check that $m_{E\circ F}\in\dcl(m_{E_i\circ F}:i\in I)$. Let $f$ be an automorphism fixing each $m_{E_i \circ F}$. Then for each $i\in I$ there is some $a_i$ such that
$$E_i(m,a_i)\wedge F(a_i,f(m)).$$
By compactness there is some $a$ such that $E(m,a)\wedge F(a,f(m))$. Hence $f$ fixes $m_{E\circ F}$.
\end{proof}

\begin{remark}  The Galois correspondence provides another proof of Corollary~\ref{LP} in terms of groups. Let $d$ be a bounded hyperimaginary and let $(e_i:i\in I)$ be a family  of finitary  normal hyperimaginaries such that $(e_i:i\in I)\sim \bdd(\e)$. 
As above, we may assume that the family is closed under finite composition. Let $H_i=\Fix(e_i)$, a closed normal subgroup of the Galois group of $T$. Under the Galois correspondence, the conditions on the $e_i$'s means that $\bigcap_i H_i=\{1\}$, and for each $i,j$ there is some $k$ such that $H_i\cap H_j=H_k$.
Let $H=\Fix (d)$, and consider $L_i=H.H_i$, a closed subgroup of the Galois group. Again, the Galois correspondence tells us that $L_i=\Fix (h_i)$ for some bounded hyperimaginary $h_i$, and certainly $h_i$ is finitary since $e_i$ is. Now $\bigcap_i L_i=\bigcap_i H.H_i=H.\bigcap_i H_i=H$, which means that $d\sim (h_i : i\in I)$.
\end{remark}

\section{Local types of hyperimaginaries}

\begin{defi}\rm Let $e,d$ be hyperimaginaries.    The \emph{orbit} of $e$ over $d$ is the set $\O(e/d)$  of all  hyperimaginaries $e^\prime$  such that $e\equiv_d e^\prime$. 
\end{defi}  

\begin{remark}  Notice that for an automorphism $f$,  the condition $f(\O(e/d))= \O(e/d)$  is equivalent the conjunction of $e\equiv_d f(e)$ and $e\equiv_d f^{-1}(e)$.
\end{remark}

Next lemma  is due to Buechler, Pillay and Wagner (Lemma 2.18 in~\cite{BuechlerPillayWagner00}). It basically says that we can consider $\O(e/d)$ as a hyperimaginary  if $e\in\bdd(d)$. In our Proposition~\ref{l4} below we have generalized this fact to any closed set in a Kim-Pillay space.  We apply this to some closed sets $\O_\varphi(e/d)$ obtaining thus some hyperimaginaries $h_{p,\varphi,d}$.  For $d\in \bdd(\emptyset)$  and $p(x)=\tp(e/\emptyset)$  we understand $\tp(e/h_{p,\varphi,d})$ as an approximation to the $\varphi$-type of $e$ over $d$. 

\begin{remark}\label{l1}  If $e\in\bdd(d)$, then $\O(e/d)$ is $\sw$-equivalent to some hyperimaginary $h$, in the sense that the automorphisms of the monster model fixing $h$ is the set of automorphisms fixing set-wise $\O(e/d)$. 
\end{remark}

\begin{lemma}\label{l2} If  $e=(a_i:i<\omega)_\KP$  and   $e_n=(a_i: i\leq n)_\KP$,  then  $e\sim (e_n: n<\omega)$.
\end{lemma}
\begin{proof} For every automorphism $f$,  $f\in \Fix(e)$   iff   $(a_i: i<\omega)\equiv_{\bdd(\emptyset)}(f(a_i):i<\omega)$  iff    $(a_i: i\leq n)\equiv_{\bdd(\emptyset)}(f(a_i):i\leq n)$ for all $n<\omega$ iff  $f(e_n)=e_n$    for all $n<\omega$.
\end{proof}

\begin{prop}\label{l3} If $d\in\bdd(\emptyset)$, then  $d\sw (\O(e/d): e\in \c^\omega/\KP)$  and  $d\sw(\O(e/d):e\in \c^n/\KP,  n<\omega )$
\end{prop}
\begin{proof}  If $f\in\Fix(d)$, then $f$  permutes every orbit $\O(e/d)$.   

Assume  $f$ permutes every $\O(e/d)$  for every countable $\KP$-class  $e\in\c^\omega/\KP$.  It is well known that each hyperimaginary is equivalent to a sequence of countable hyperimaginaries. Hence $d\sim (d_i:i\in I)$, where every $d_i$ is a countable hyperimaginary. Choose an $\omega$-tuple $a_i$ and a bounded 0-type-definable equivalence relation $E_i$ such that  $d_i={a_i}_{E_i}$.   By hypothesis,  $f({a_i}_\KP)\in \O({a_i}_\KP/d)$ and therefore  $f({a_i}_\KP)=g_i({a_i}_\KP)$ for some $g_i\in \Fix(d)$.   Note that $d_i$ is a union of $\KP$-classes of $\omega$-tuples.  Since $g_i$ fixes $d_i$,  $f$ permutes these $\KP$-classes and then $f(d_i)=d_i$.  Since $f$ fixes each $d_i$,  $f(d)=d$.

Assume now $f$ permutes every $\O(e/d)$ for every finitary  $\KP$-class  $e\in \c^n/\KP$.  We show that $f$ permutes $\O(e/d)$ for every countable $\KP$-class  $e\in \c^\omega/\KP$.  Let $e=(a_i:i<\omega)_\KP$ and let  $e_n= (a_i:i\leq n)_\KP$.   Since  $e_n\equiv_d  f(e_n)$  for all $n<\omega$,   $(e_n: n<\omega)\equiv_d (f(e_n): n<\omega)$.  Choose $g\in\Fix(d)$ such that $g(e_n:n<\omega)=(f(e_n): n<\omega)$.  Then  $f^{-1}g(e_n)=e_n$  for all $n<\omega$  and by Lemma~\ref{l2}  $f^{-1}g(e)=e$. It follows that $e\equiv_d f(e)$ and hence $f$ permutes $\O(e/d)$.
\end{proof}

\begin{prop}\label{l4} Every closed set $C$ in a Kim-Pillay space is $\sw$-equivalent to a hyperimaginary $h_C$, that is, the automorphisms of the monster model fixing set-wise $C$ are the automorphisms fixing $h_C$.
\end{prop}
\begin{proof}  Let  $E$ be a bounded 0-type-definable equivalence relation on $\alpha$-tuples and let $X=\c^\alpha/E$ be the corresponding Kim-Pillay space. If $C\subseteq X$ is closed, then for some partial type $\pi(x,z)$, for some tuple $b$, $\pi(\c,b)=\{a: a_E\in C\}$. For each formula $\theta(x,y)\in E(x,y)$ there is a maximal length $n=n_\theta<\omega$ of a sequence of tuples  $(a_i:i<n)$ such that ${a_i}_E\in C$ and  $\models \neg \theta(a_i,a_j)$  for all $i<j<n$.  Let $(a^\theta_i:i<n_\theta)$ witness it,  let  $\Sigma_\theta(z,z^\prime)$ be the partial type
$$\exists (x_i: i<n_\theta) ( \bigwedge_{i<j<{n_\theta}}\neg \theta(x_i,x_j)\wedge \bigwedge_{i<n_\theta}\pi(x_i,z)\wedge \bigwedge_{i<n_\theta}\pi(x_i,z^\prime))$$
and 
$$F(z,z^\prime) =  \bigwedge_{\theta\in E} \Sigma_\theta(z,z^\prime)$$
\emph{Claim}: For every automorphism $f$,  $f(C)=C$  if and only if   $\models F(b,f(b))$.\\
\emph{Proof of the claim}:  From left to right it is straightforward. For the other direction, assume  $\models F(b,f(b))$ and choose  $(c^\theta_i: i<n_\theta,\theta\in E)$ witnessing it.   Let $a_E\in C$.  Then  $\models \pi(a,b)$, and hence  $\models \pi(f(a),f(b))$. By maximality of $n_\theta$, for every $\theta\in E$ there is some $i<n_\theta$ such that  $\models\theta(f(a), c^\theta_i)$. By compactness,   $E(f(a),x) \cup \pi(x,b)$ is consistent and therefore  $f(a_E)\in C$. This shows that  $f(C)\subseteq C$.   By the same reason, $f^{-1}(C)\subseteq  C$, that is, $C\subseteq f(C)$.

It follows from the claim that $F$ defines a 0-type-definable equivalence relation on realizations of $p(x)=\tp(b)$. By standar arguments it can be extended to a 0-type-definable equivalence relation defined for all tuples of the length of $b$.  The hyperimaginary $b_F$ satisfies the requirements.
\end{proof}

\begin{defi}\rm  Let $e,d$ be hyperimaginaries. If $\varphi(x,y)\in L$,  $\models \varphi(e,d)$  means that  $\models \varphi(a,b)$  for some representatives $a,b$ of $e,d$ respectively.   Notice that $e\equiv_d e^\prime$ iff  $\models \varphi(e,d)\Leftrightarrow \models \varphi(e^\prime,d)$  for all $\varphi(x,y)\in L$. 
 Let $\O_\varphi(e/d)=\{e^\prime:  e^\prime\equiv e \mbox{ and} \models \varphi(e^\prime,d)\}$.  Let $p(x)=\tp(e)$ and assume $e\in\bdd(\emptyset)$.   Then $e=a_E$  for some tuple $a$ and some bounded 0-type-definable equivalence relation $E$.  The set of all $E$-classes is a Kim-Pillay space and $\O_\varphi(e/d)$ defines a closed subset. By Lemma~\ref{l4}  there is some hyperimaginary $h_{p,\varphi,d}$  such that $$h_{p,\varphi,d}\sw\O_\varphi(e/d).$$
\end{defi}

 The equivalence relation $E(e,e^\prime)$ defined by $ \models \varphi(e,d) \Leftrightarrow \models \varphi(e^\prime,d)$ is not, in general, type-definable. This is the reason why an adequate treatment of local types (or $\varphi$-types) is missing in the model theory of hyperimaginaries. The following results show that the types  $\tp(e/h_{\tp(e),\varphi,d})$  are  (for $e$ bounded) a substitute for the $\varphi$-type of $e$ over $d$  and we apply this in Corollary~\ref{l7}  to obtain a new decomposition of a bounded hyperimaginary in terms of orbits.

\begin{remark}\label{l5} Let $d$ be a hyperimaginary, $e\in \bdd(\emptyset)$,  $p(x)=\tp(e)$ and $\varphi(x,y)\in L$.
\begin{enumerate}
\item If $e^\prime\equiv_{h_{p,\varphi,d}}e$ then  $\models \varphi(e,d)\Leftrightarrow \models \varphi(e^\prime,d)$.
\item $h_{p,\varphi,d}\in\dcl(d)$.
\end{enumerate}
\end{remark}
\begin{proof}  Clear.
\end{proof}

\begin{prop}\label{l6} Let $d$ be a hyperimaginary, $e\in \bdd(\emptyset)$  and $p(x)=\tp(e)$.  For any $e^\prime\models p$:
$$e^\prime\equiv_d e   \mbox{ if and only if } e^\prime\equiv_{h_{p,\varphi,d}}e \mbox{ for every } \varphi(x,y)\in L$$
\end{prop}
\begin{proof}
By Remark~\ref{l5}.
\end{proof}  

\begin{cor}\label{l7} If $d\in\bdd(\emptyset)$, then $d\sw(\O(e/h_{\tp(e),\varphi,d}):  e\in\c^n/\KP, n<\omega)$.
\end{cor}
\begin{proof} Let $d\in \bdd(\emptyset)$.  If $f\in \Fix(d)$, then $f$ fixes  $h_{\tp(e),\varphi,d}$ and permutes  $\O(e/h_{\tp(e),\varphi,d})$. On the other hand, if $e\in \c^n/\KP$ and $f$ permutes all the orbits $(\O(e/h_{\tp(e),\varphi,d})$, then  $e\equiv_{h_{\tp(e),\varphi,d}}f(e)$ for all $\varphi$  and by Proposition~\ref{l6}  $e\equiv_d f(e)$.  Similarly, $e\equiv_d f^{-1}(e)$. It follows  that  $f$ permutes $\O(e/d)$.   By
Proposition~\ref{l3}, $f(d)=d$.
\end{proof}


\noindent
 Departamento de L\'ogica, Historia y Filosof\'{\i}a de la  Ciencia\\
 Universidad de Barcelona\\
 Montalegre 6, 08001 Barcelona\\
 ESPA\~NA\\
 E-mail: {\tt e.casanovas@ub.edu}  and  {\tt z.-@laposte.net}


\begin{thebibliography}{1}

\bibitem{BuechlerPillayWagner00}
S.~Buechler, A.~Pillay, and F.~O. Wagner.
\newblock Supersimple theories.
\newblock {\em Journal of the American Mathematical Society}, 14:109--124,
  2000.

\bibitem{Casanovas09}
E.~Casanovas.
\newblock {\em Simple theories and hyperimaginaries}, volume~39 of {\em Lecture
  Notes in Logic}.
\newblock Cambridge University Press, 2011.

\bibitem{Las-Pil98}
D.~Lascar and A.~Pillay.
\newblock {H}yperimaginaries and automorphism groups.
\newblock {\em The Journal of Symbolic Logic}, 66:127--143, 2001.

\bibitem{Ziegler01}
M.~Ziegler.
\newblock Introduction to the {L}ascar group.
\newblock In K.~Tent, editor, {\em Tits Buildings and the Model Theory of
  Groups}, number 291 in London Mathematical Society Lecture Notes Series,
  pages 279--298. Cambridge University Press, Cambridge, 2002.

\end{thebibliography}
\end{document}